\documentclass{amsart}

\usepackage[leqno]{amsmath}
\usepackage{amsthm,latexsym,amsfonts,amssymb}
\usepackage[all]{xy} \SelectTips{eu}{}
\usepackage{hyperref}

\newcommand{\numberseries}{\bfseries}   
\newlength{\thmtopspace}                
\newlength{\thmbotspace}                
\newlength{\thmheadspace}               
\newlength{\thmindent}                  
\newtheoremstyle{bfupright head,slanted body}
                {\thmtopspace}{\thmbotspace}
                {\slshape}{\thmindent}{\bfseries}{.}{\thmheadspace}
                {{\numberseries \thmnumber{#2\;}}\thmnote{#3}}

\newtheoremstyle{fixed bf head,slanted body}
                {\thmtopspace}{\thmbotspace}{\slshape}
                {\thmindent}{\bfseries}{.}{\thmheadspace}
                {{\numberseries \thmnumber{#2\;}}\thmname{#1}\thmnote{ (#3)}}

\newtheoremstyle{fixed bf head,upright body}
                {\thmtopspace}{\thmbotspace}{\upshape}
                {\thmindent}{\bfseries}{.}{\thmheadspace}
                {{\numberseries \thmnumber{#2\;}}\thmname{#1}\thmnote{ (#3)}}

\newtheoremstyle{numbered paragraph}
                {\thmtopspace}{\thmbotspace}{\upshape}
                {\thmindent}{\upshape}{}{\thmheadspace}
                {{\numberseries \thmnumber{#2.}}}

\theoremstyle{bfupright head,slanted body}
\newtheorem{res}{}[section]
\newtheorem*{res*}{}
\theoremstyle{fixed bf head,slanted body}
\newtheorem{thm}[res]{Theorem}          \newtheorem*{thm*}{Theorem}
\newtheorem{prp}[res]{Proposition}      \newtheorem*{prp*}{Proposition}
\newtheorem{cor}[res]{Corollary}        \newtheorem*{cor*}{Corollary}
\newtheorem{lem}[res]{Lemma}            \newtheorem*{lem*}{Lemma}

\theoremstyle{fixed bf head,upright body}
\newtheorem{dfn}[res]{Definition}       \newtheorem*{dfn*}{Definition}
\newtheorem{con}[res]{Construction}     \newtheorem*{con*}{Construction}
\newtheorem{rmk}[res]{Remark}           \newtheorem*{rmk*}{Remark}
\newtheorem{fct}[res]{Fact}             \newtheorem*{fct*}{Fact}

\theoremstyle{numbered paragraph}
\newtheorem{ipg}[res]{}

\newlength{\thmlistleft}        
\newlength{\thmlistright}       
\newlength{\thmlistpartopsep}   
\newlength{\thmlisttopsep}      
\newlength{\thmlistparsep}      
\newlength{\thmlistitemsep}     

\newcounter{eqc}
\newenvironment{eqc}{\begin{list}{\upshape (\textit{\roman{eqc}})}%
    {\usecounter{eqc}%
      \setlength{\leftmargin}{\thmlistleft}%
      \setlength{\labelwidth}{\thmlistleft}%
      \setlength{\rightmargin}{\thmlistright}%
      \setlength{\partopsep}{\thmlistpartopsep}%
      \setlength{\topsep}{\thmlisttopsep}%
      \setlength{\parsep}{\thmlistparsep}%
      \setlength{\itemsep}{\thmlistitemsep}}}%
  {\end{list}}%
\newcommand{\eqclbl}[1]{{\upshape(\textit{#1})}}

\newcounter{prt}
\newenvironment{prt}{\begin{list}{\upshape (\alph{prt})}%
    {\usecounter{prt}%
      \setlength{\leftmargin}{\thmlistleft}%
      \setlength{\labelwidth}{\thmlistleft}%
      \setlength{\rightmargin}{\thmlistright}%
      \setlength{\partopsep}{\thmlistpartopsep}%
      \setlength{\topsep}{\thmlisttopsep}%
      \setlength{\parsep}{\thmlistparsep}%
      \setlength{\itemsep}{\thmlistitemsep}}}%
  {\end{list}}%

\newcounter{rqm}
\newenvironment{rqm}{\begin{list}{\upshape (\arabic{rqm})}%
    {\usecounter{rqm}%
      \setlength{\leftmargin}{\thmlistleft}%
      \setlength{\labelwidth}{\thmlistleft}%
      \setlength{\rightmargin}{\thmlistright}%
      \setlength{\partopsep}{\thmlistpartopsep}%
      \setlength{\topsep}{\thmlisttopsep}%
      \setlength{\parsep}{\thmlistparsep}%
      \setlength{\itemsep}{\thmlistitemsep}}}%
  {\end{list}}%

\newenvironment{prf*}[1][Proof]{%
  \begin{proof}[\bf #1]
    \setcounter{equation}{0}
    }
  {\end{proof}
}
   \newcommand{\proofofimp}[3][:]{\mbox{\eqclbl{#2}$\!\implies\!$\eqclbl{#3}#1}}

\newcommand{\pgref}[1]{\ref{#1}}
\newcommand{\fctref}[2][Fact~]{#1\ref{fct:#2}}
\newcommand{\thmref}[2][Theorem~]{#1\pgref{thm:#2}}
\newcommand{\corref}[2][Corollary~]{#1\pgref{cor:#2}}
\newcommand{\prpref}[2][Proposition~]{#1\pgref{prp:#2}}
\newcommand{\lemref}[2][Lemma~]{#1\pgref{lem:#2}}

\newcommand{\dfnref}[2][Definition~]{#1\pgref{dfn:#2}}
\newcommand{\rmkref}[2][Remark~]{#1\pgref{rmk:#2}}
\renewcommand{\eqref}[1]{(\pgref{eq:#1})}
\newcommand{\thmcite}[2][?]{\cite[Thm.~#1]{#2}}
\newcommand{\rmkcite}[2][?]{\cite[Rmk.~#1]{#2}}
\newcommand{\corcite}[2][?]{\cite[Cor.~#1]{#2}}
\newcommand{\prpcite}[2][?]{\cite[Prop.~#1]{#2}}
\newcommand{\lemcite}[2][?]{\cite[Lem.~#1]{#2}}
\newcommand{\seccite}[2][?]{\cite[Sect.~#1]{#2}}
\newcommand{\dfncite}[2][?]{\cite[Def.~#1]{#2}}
\newcommand{\exacite}[2][?]{\cite[Exa.~#1]{#2}}
\def\urltilda{\kern -.15em\lower .7ex\hbox{\~{}}\kern .04em}
\newcommand{\setof}[3][\mspace{1mu}]{\{#1#2 \mid #3#1\}}
\newcommand{\ZZ}{\mathbb{Z}}
\newcommand{\qtext}[1]{\quad\text{#1}\quad}
\newcommand{\qqtext}[1]{\qquad\text{#1}\qquad}
\newcommand{\qand}{\qtext{and}}
\newcommand{\qqand}{\qqtext{and}}
\newcommand{\deq}{\:=\:}
\newcommand{\dle}{\:\le\:}
\newcommand{\dge}{\:\ge\:}
\newcommand{\gra}{\alpha}
\newcommand{\grb}{\beta}
\newcommand{\grg}{\gamma}
\newcommand{\gre}{\varepsilon}
\newcommand{\is}{\cong}
\newcommand{\qis}{\simeq}
\renewcommand{\le}{\leqslant}
\renewcommand{\ge}{\geqslant}

\newcommand{\lra}{\longrightarrow}
\newcommand{\xla}[2][]{\xleftarrow[#1]{\;#2\;}}
\newcommand{\xra}[2][]{\xrightarrow[#1]{\;#2\;}}
\newcommand{\qla}{\xla{\qis}}
\newcommand{\qra}{\xra{\qis}}
\newcommand{\QQ}{\mathbb{Q}}
\newcommand{\QZ}{\QQ/\ZZ}
\newcommand{\Rop}{R^\circ}
\newcommand{\mapdef}[4][\rightarrow]{\nobreak{#2\colon #3 #1 #4}}
\newcommand{\qisdef}[4][\xra{\qis}]{\nobreak{#2\colon #3 #1 #4}}

\renewcommand{\Im}[1]{\nobreak{\operatorname{Im}#1}}
\newcommand{\Ker}[1]{\nobreak{\operatorname{Ker}#1}}
\newcommand{\Coker}[1]{\nobreak{\operatorname{Coker}#1}}
\newcommand{\Cone}[1]{\nobreak{\operatorname{Cone}#1}}
\newcommand{\cls}[1]{[#1]}
\newcommand{\dif}[2][]{{\partial}^{#2}_{#1}}
\newcommand{\Bo}[2][]{\operatorname{B}_{#1}(#2)}
\newcommand{\Cy}[2][]{\operatorname{Z}_{#1}(#2)}
\newcommand{\Co}[2][]{\operatorname{C}_{#1}(#2)}
\renewcommand{\H}[2][]{\operatorname{H}_{#1}(#2)}
\newcommand{\Shift}[2][]{\mathsf{\Sigma}^{#1}{#2}}
\newcommand{\Tha}[2]{#2_{{\scriptscriptstyle\le}#1}}
\newcommand{\Thb}[2]{#2_{{\scriptscriptstyle\ge}#1}}
\newcommand{\Tsa}[2]{#2_{{\scriptscriptstyle\subseteq}#1}}
\newcommand{\Tsb}[2]{#2_{{\scriptscriptstyle\supseteq}#1}}
\newcommand{\fd}[2][R]{\operatorname{fd}_{#1}#2}
\newcommand{\id}[2][R]{\operatorname{id}_{#1}#2}
\newcommand{\Gfd}[2][R]{\operatorname{Gfd}_{#1}#2}
\newcommand{\Gfcd}[2][R]{\operatorname{Gfcd}_{#1}#2}
\newcommand{\Hom}[3][R]{\operatorname{Hom}_{#1}(#2,#3)}
\newcommand{\RHom}[3][R]{\operatorname{\mathbf{R}Hom}_{#1}(#2,#3)}
\newcommand{\Ext}[4][R]{\operatorname{Ext}_{#1}^{#2}(#3,#4)}

\newcommand{\tp}[3][R]{\nobreak{#2\otimes_{#1}#3}}

\makeatletter
\def\@nobreak@#1{\mathchoice%
  {\nobreakdef@\displaystyle\f@size{#1}}%
  {\nobreakdef@\nobreakstyle\tf@size{\firstchoice@false #1}}%
  {\nobreakdef@\nobreakstyle\sf@size{\firstchoice@false #1}}%
  {\nobreakdef@\nobreakstyle\ssf@size{\firstchoice@false #1}}%
  \check@mathfonts}%
\def\nobreakdef@#1#2#3{\hbox{{%
                    \everymath{#1}%
                    \let\f@size#2\selectfont%
                    #3}}}%
\makeatother

\newcommand{\Ftac}{\textnormal{\bf F}-totally acyclic}

\newcommand{\Text}[4][{\sf FC}]{\smash{\operatorname{\widehat{Ext}}%
  }_{#1}^{#2^{\phantom{|}\mspace{-6mu}}}(#3,#4)}
\newcommand{\dc}[2]{\operatorname{D}^{#1}(#2)}

\numberwithin{equation}{res}

\begin{document}

\title[A refinement of Gorenstein flat dimension]{%
  A refinement of Gorenstein flat dimension\\via the flat--cotorsion
  theory}

\author[L.W.\ Christensen]{Lars Winther Christensen}

\address{L.W.C. \ Texas Tech University, Lubbock, TX 79409, U.S.A.}

\email{lars.w.christensen@ttu.edu}

\urladdr{http://www.math.ttu.edu/\urltilda lchriste}

\author[S.\ Estrada]{Sergio Estrada}

\address{S.E. \ Universidad de Murcia, Murcia 30100, Spain}

\email{sestrada@um.es}

\urladdr{http://webs.um.es/sestrada}

\author[L.\ Liang]{Li Liang}

\address{L.L. \ Lanzhou Jiaotong University, Lanzhou 730070, China}

\email{lliangnju@gmail.com}

\urladdr{https://sites.google.com/site/lliangnju}

\author[P.\ Thompson]{Peder Thompson}

\address{P.T. \ Norwegian University of Science and Technology, 7491
  Trondheim, Norway}

\email{peder.thompson@ntnu.no}

\urladdr{https://folk.ntnu.no/pedertho}

\author[D.\ Wu]{Dejun Wu}

\address{D.W. \ Lanzhou University of Technology, Lanzhou 730050,
  China}

\email{wudj@lut.cn}

\author[G.\ Yang]{Gang Yang}

\address{G.Y. \ Lanzhou Jiaotong University, Lanzhou 730070, China}

\email{yanggang@mail.lzjtu.cn}

\thanks{L.W.C.\ was partly supported by Simons Foundation
  collaboration grant 428308; S.E.\ was partly supported by grant
  MTM2016-77445-P (AEI/FEDER,UE) and Fundaci\'on Seneca grant
  19880/GERM/15; L.L.\ was partly supported by NSF of China grant
  11761045 and NSF of Gansu Province grant 18JR3RA113; D.W.\ was
  partly supported by NSF of China grants 11761047 and 11861043; G.Y\
  was partly supported by NSF of China grant 11561039; both L.L\ and
  G.Y\ were also partly supported by the Foundation of A Hundred Youth
  Talents Training Program of Lanzhou Jiaotong University. The
  research was partly done during D.W.'s year-long visit to Texas Tech
  University; the hospitality of the TTU Department of Mathematics and
  Statistics is acknowledged with gratitude.}

\date{28 September 2020}



\begin{abstract}
  We introduce a refinement of the Gorenstein flat dimension for
  complexes over an associative ring---the \emph{Gorenstein
    flat-cotorsion dimension}---and prove that it, unlike the
  Gorenstein flat dimension, behaves as one expects of a homological
  dimension without extra assumptions on the ring. Crucially, we show
  that it coincides with the Gorenstein flat dimension for complexes
  where the latter is finite, and for complexes over right coherent
  rings---the setting where the Gorenstein flat dimension is known to
  behave as expected.
\end{abstract}

\maketitle

\thispagestyle{empty}

\section*{Introduction}

\noindent
The introduction of the G-dimension by Auslander and Bridger
\cite{MAsMBr69}, and the subsequent broader notion of Gorenstein
projective dimension by Enochs and Jenda \cite{EEnOJn95b}, provided
for an elegant characterization of Gorenstein rings in terms of
finiteness of homological invariants. It is modeled on the
characterization by Auslander, Buchsbaum, and Serre of commutative
regular local rings as rings of finite global dimension. Further
pursuit of this analogy led to the introduction of the Gorenstein
injective and Gorenstein flat dimensions with the aim of building a
theory of Gorenstein homological dimensions modeled on the classic
projective, injective, and flat dimensions.

This program has largely been successful, but from a homological
algebra point of view not entirely so: While the Gorenstein projective
and injective dimensions behave as one expects of homological
dimensions---in particular, they can be computed in terms of vanishing
of cohomology---the Gorenstein flat dimension only exhibits such
behavior under coherence assumptions on the ring; see Holm
\cite{HHl04a}. Our goal is to give a new perspective on the Gorenstein
flat dimension: One that puts it on the same footing as the Gorenstein
projective and Gorenstein injective dimensions and ensures that it
behaves nicely without assumptions on the ring.

In recent years, it has become apparent that one ought to pay special
attention to the narrower class of Gorenstein flat modules that are
also cotorsion. For example, work of Gillespie \cite{JGl17} shows that
under coherence assumptions on the ring, the category of modules that
are Gorenstein flat and cotorsion is Frobenius, while the category of
Gorenstein flat modules rarely is; in fact, it only happens when every
Gorenstein flat module \textsl{is} cotorsion, see
\thmcite[4.5]{CETa}. Motivated in part by this, \emph{Gorenstein
  flat-cotorsion} modules were introduced in \cite{CETa}, and it was
shown that over a right coherent ring they are precisely the modules
that are Gorenstein flat and cotorsion.

To push the Gorenstein flat dimension beyond the setting of coherent
rings, we introduce the \emph{Gorenstein flat-cotorsion dimension,}
not so much to introduce a new homological dimension but rather to
refine the already established dimension. The Gorenstein
flat-cotorsion dimension is defined in terms of the Hom
functor---rather than the tensor product functor which is used for the
Gorenstein flat dimension---and in this way it behaves more like the
Gorenstein projective and Gorenstein injective dimensions.

Let $R$ be an associative ring. The flat, Gorenstein flat, and
Gorenstein flat-cotorsion dimensions of an $R$-complex $M$ are written
$\fd{M}$, $\Gfd{M}$ and $\Gfcd{M}$, respectively. The next two
statements, which are extracted from \thmref[Theorems~]{Gfcd},
\thmref[]{Gfd-Gfcd}, and \thmref[]{fd-Gfcd}, capture the essence of
the new dimension.

\begin{res*}[Theorem A]
  \label{res:Gfcd}
  Let $M$ be an $R$-complex and $n$ an integer. If $\Gfcd{M}$ is
  finite and $n \ge \sup{M}$, then the following conditions are
  equivalent.
  \begin{eqc}
  \item $\Gfcd{M} \le n$.
  \item $\Ext{i}{M}{C} =0$ for all $i > n$ and every $R$-module $C$
    that is flat and cotorsion.
  \item $\Ext{i}{M}{C} =0$ for all $i > n$ and every cotorsion
    $R$-module $C$ of finite flat dimension.
  \item $\Ext{n+1}{M}{C} =0$ for every cotorsion $R$-module $C$ of
    finite flat dimension.
  \end{eqc}
\end{res*}

\begin{res*}[Theorem B]
  \label{res:Gfd-Gfcd}
  Let $M$ be an $R$-complex. There are inequalities,
  \begin{equation*}
    \Gfcd{M} \dle \Gfd{M} \dle \fd{M}\:,
  \end{equation*}
  and if any of these quantities is finite, then it equals those to
  the left of it.
\end{res*}

While the proof of Theorem~A is standard fare homological algebra, the
connection to the Gorenstein flat dimension captured by Theorem B
relies crucially on recent work of \v{S}aroch and
\v{S}tov{\'{\i}}{\v{c}}ek \cite{JSrJSt20}. It is their work that
allows us to say that the Gorenstein flat-cotorsion dimension is not
so much a new homological dimension as it is a new perspective on an
old one.  To further illustrate the utility of this perspective, we
briefly introduce a version of Tate cohomology associated to the
theory of Gorenstein flat-cotorsion modules. It is shown that this
yields a generalization of a result of Hu and Ding \cite{JHuNDn16}
that characterizes complexes of finite flat dimension among those of
finite Gorenstein flat dimension; see \thmref{vanishing}.

\section{Notation and Terminology}

\noindent
Throughout the paper, $R$ denotes an associative ring. By an
$R$-module we mean a left $R$-module; right $R$-modules are considered
modules over the opposite ring $\Rop$.  A complex of $R$-modules is,
for short, called an $R$-complex. We use homological notation for
complexes, i.e.\ for $n\in\ZZ$ the module in degree $n$ of an
$R$-complex $M$ is denoted $M_n$. The submodules of boundaries and
cycles are denoted $\Bo[n]{M}$ and $\Cy[n]{M}$, respectively. The
homology is as always the quotient $\H[n]{M} =
\Cy[n]{M}/\Bo[n]{M}$. We further use the notation $\Co[n]{M}$ for the
cokernel of the differential $\dif[n+1]{M}$, i.e.\ one has
$\Co[n]{M} = M_n/\Bo[n]{M}$.  The invariants
\begin{equation*}
  \sup{M} = \sup\setof{n\in\ZZ}{\H[n]{M} \ne 0} \qand
  \inf{M} = \inf\setof{n\in\ZZ}{\H[n]{M} \ne 0}
\end{equation*}
capture the homological position of the complex. If $\H[n]{M}=0$ holds
for all $n\in\ZZ$, then $M$ is called \emph{acyclic.} Morphisms of
complexes that induce isomorphisms in homology are called
\emph{quasi-isomorphisms;} they are characterized by having acyclic
mapping cones. For $s\in\ZZ$ the $s$-fold \emph{shift} of $M$ is the
complex $\Shift[s]{M}$ defined by $(\Shift[s]{M})_n = M_{n-s}$ and
$\dif[n]{\Shift[s]{M}} = (-1)^s\dif[n-s]{M}$.

For an $R$-complex $M$, the \emph{hard truncation above} of $M$ at $n$
is the complex
\begin{equation*}
  \Tha{n}{M} = \cdots \lra 0 \lra M_n \lra M_{n-1} \lra  M_{n-2} \lra \cdots
\end{equation*}
with differential induced from $M$; the \emph{hard truncation below}
of $M$ at $n$, denoted $\Thb{n}{M}$, is defined similarly. The
\emph{soft truncation above} of $M$ at $n$ is the complex
\begin{equation*}
  \Tsa{n}{M}= \cdots \lra 0 \lra \Co[n]{M} \lra M_{n-1} \lra M_{n-2}\lra\cdots
\end{equation*}
with differential induced from $M$.  The \emph{soft truncation below}
of $M$ at $n$, denoted $\Tsb{n}{M}$, is defined similarly by replacing
$M_n$ by $\Cy[n]{M}$.

An $R$-complex $P$ is called \emph{semi-projective} if it consists of
projective modules and the functor $\Hom{P}{-}$ preserves
acyclicity. Dually, a complex $I$ is called \emph{semi-injective} if
it consists of injective modules and the functor $\Hom{-}{I}$
preserves acyclicity. Every $R$-complex $M$ has a semi-projective
resolution, i.e.\ there is a quasi-isomorphism $\pi\colon P \qra M$
where $P$ is semi-projective; one can choose $\pi$ surjective or one
can choose $P$ with $P_n=0$ for $n < \inf{M}$; see Avramov and Foxby
\seccite[1]{LLAHBF91}. Dually, every $R$-complex $M$ has a
semi-injective resolution i.e.\ there is a quasi-isomorphism
$\iota\colon M \qra I$ where $I$ is semi-injective; one can choose
$\iota$ injective or one can choose $I$ with $I_n=0$ for
$n > \sup{M}$.

An $R$-complex $F$ is called \emph{semi-flat} if it consists of flat
modules and the functor $\tp{-}{F}$ preserves acyclicity; see
\seccite[1]{LLAHBF91}. Every semi-projective complex is semi-flat and
so is every bounded below complex of flat modules.

An $R$-complex $C$ is called \emph{semi-cotorsion} if it consists of
cotorsion modules and $\Hom{F}{C}$ is acyclic for every acyclic
semi-flat $R$-complex $F$. Every semi-injective complex is
semi-cotorsion and so is every bounded above complex of cotorsion
$R$-modules; this is a standard consequence of
\lemcite[2.5]{CFH-06}. In fact, one can do better:

\begin{fct}
  \label{fct:BCE}
  An acyclic complex of cotorsion $R$-modules has cotorsion cycle
  modules, and it follows that every complex of cotorsion $R$-modules
  is semi-cotorsion. This is proved by Bazzoni, Cort\'es Izurdiaga,
  and Estrada \thmcite[1.3]{BCE}.
\end{fct}

\begin{lem}
  \label{lem:cot}
  For every $\Rop$-complex $M$ the $R$-complex $\Hom[\ZZ]{M}{\QZ}$ is
  semi-cotorsion.
\end{lem}

\begin{prf*}
  Let $F$ be an acyclic semi-flat $R$-complex; it is a direct limit of
  contractible complexes of finitely generated free $R$-modules, see
  \thmcite[7.3]{LWCHHl15}, so $\tp{M}{F}$ is acyclic. It follows by
  Hom--tensor adjunction that $\Hom[\ZZ]{M}{\QZ}$ is a complex of
  cotorsion $R$-modules and, further, that there is an isomorphism
  \begin{equation*}
    \Hom{F}{\Hom[\ZZ]{M}{\QZ}} \is \Hom[\ZZ]{\tp{M}{F}}{\QZ}\:.
  \end{equation*}
  The right-hand complex is acyclic, so $\Hom[\ZZ]{M}{\QZ}$ is
  semi-cotorsion.
\end{prf*}

\begin{lem}
  \label{lem:3-cot}
  Let $0 \to C' \to C \to C'' \to 0$ be an exact sequence of
  $R$-complexes. If $C'$ is semi-cotorsion, then $C$ is semi-cotorsion
  if and only if $C''$ is semi-cotorsion.
\end{lem}

\begin{prf*}
  This is an immediate consequence of \fctref{BCE}.
\end{prf*}

The next fact is standard and proved in much the same way as
\lemref{3-cot}.

\begin{fct}
  \label{fct:3-sf}
  Let $0 \to F' \to F \to F'' \to 0$ be an exact sequence of
  $R$-complexes. If $F''$ is semi-flat, then $F$ is semi-flat if and
  only if $F'$ is semi-flat.
\end{fct}

From \cite[A.1]{TNkPTh} we recall:

\begin{fct}
  \label{fct:rep}
  For a semi-flat $R$-complex $F$ and a semi-cotorsion $R$-complex $C$
  there is an isomorphism $\RHom{F}{C} \qis \Hom{F}{C}$ in the derived
  category.
\end{fct}

\begin{prp}
  \label{prp:semi-flat}
  A complex $F$ of flat $R$-modules is semi-flat if and only if the
  complex $\Hom{F}{C}$ is acyclic for every acyclic semi-cotorsion
  complex $C$.
\end{prp}

\begin{prf*}
  The ``only if'' follows from \fctref{rep}. For the converse, let $M$
  be an acyclic $\Rop$-complex. The $R$-complex $\Hom[\ZZ]{M}{\QZ}$ is
  semi-cotorsion by \lemref{cot} and acyclic. Thus it follows from the
  adjunction isomorphism
  $$\Hom[\ZZ]{\tp{M}{F}}{\QZ} \is \Hom{F}{\Hom[\ZZ]{M}{\QZ}}$$ and
  faithful injectivity of $\QZ$ that $\tp{M}{F}$ is acyclic.
\end{prf*}

\begin{prp}
  \label{prp:cothty}
  Let $C$ be a semi-cotorsion $R$-complex and $\qisdef{\grb}{F}{F'}$ a
  quasi-isomorphism of semi-flat $R$-complexes. For every morphism
  $\mapdef{\gra}{F}{C}$ there is a morphism $\mapdef{\grg}{F'}{C}$
  with $\grg\grb \sim \gra$, and $\grg$ is unique up to homotopy.
\end{prp}

\begin{prf*}
  The mapping cone of $\grb$ is acyclic and semi-flat, so the induced
  morphism $\Hom{\grb}{C}$ is a quasi-isomorphism. It follows that
  there exists a morphism $\grg$ in $\Cy[0]{\Hom{F'}{C}}$ such that
  \begin{equation*}
    \cls{\gra} \deq \H[0]{\Hom{\grb}{C}}(\cls{\grg}) \deq \cls{\grg\grb}\;;
  \end{equation*}
  that is, $\gra - \grg\grb$ is in $\Bo[0]{\Hom{F}{C}}$. Given another
  morphism $\grg'$ such that \mbox{$\grg'\grb \sim \gra$}, one has
  $\cls{\gra} = \cls{\grg'\grb}$ and, therefore
  $0 = \cls{(\grg-\grg')\grb} = \H[0]{\Hom{\grb}{C}}(\cls{\grg -
    \grg'})$. It follows that the homology class $\cls{\grg-\grg'}$ is
  $0$ as $\H[0]{\Hom{\grb}{C}}$ is an isomorphism, so $\grg-\grg'$ is
  in $\Bo[0]{\Hom{F'}{C}}$. That is, $\grg$ and $\grg'$ are homotopic.
\end{prf*}

\section{Semi-flat-cotorsion complexes}

\noindent
We recall from \cite{TNkPTh} that an $R$-complex $W$ is referred to
as \emph{semi-flat-cotorsion} if it is semi-flat and semi-cotorsion.

\begin{prp}
  \label{prp:cot-hty}
  Let $W$ be a semi-flat-cotorsion $R$-complex and $F$ be a semi-flat
  $R$-complex. If $\mapdef{\grb}{W}{F}$ is a quasi-isomorphism, then
  there is a quasi-isomorphism $\mapdef{\grg}{F}{W}$ such that
  $\grg\grb \sim 1^W$. In particular, a quasi-isomorphism of
  semi-flat-cotorsion $R$-complexes is a homotopy equivalence.
\end{prp}

\begin{prf*}
  For the first assertion apply \prpref{cothty} with $\gra = 1^W$ to
  get a morphism $\mapdef{\grg}{F}{W}$ with $\grg\grb \sim 1^W$. As
  $\grb$ and $1^W$ are quasi-isomorphisms, so is $\grg$. Next notice
  that if $F$ too is semi-flat-cotorsion then the same argument
  applies to yield a morphism $\mapdef{\grb'}{W}{F}$ with
  $\grb'\grg \sim 1^{F}$. It follows that $\grg$ and hence $\grb$ is a
  homotopy equivalence.
\end{prf*}

Gillespie~\cite{JGl04} studies how a cotorsion pair in the category of
modules induces cotorsion pairs in the category of complexes. The
short exact sequences below are often referred to as
\emph{approximations;} they exist\footnote{Although
  \cite[cor. 4.10]{JGl04} is stated for commutative rings, it is
  standard that the result remains valid without this assumption; see
  for example the discussion before \thmcite[4.2]{SEsJGl19}.} by
\corcite[4.10]{JGl04}.

\begin{fct}
  \label{fct:gil}
  For every $R$-complex $M$ there are exact sequences of $R$-complexes
  \begin{equation*}
    0 \lra M \lra C \lra F \lra 0 \qqand 0 \lra C' \lra F' \lra M \lra 0
  \end{equation*}
  where $C$ and $C'$ are semi-cotorsion, $F$ and $F'$ are semi-flat,
  and $F$ and $C'$ are acyclic.
\end{fct}

\begin{dfn}
  Let $M$ be an $R$-complex. A semi-flat-cotorsion complex that is
  isomorphic to $M$ in the derived category is called a
  \emph{semi-flat-cotorsion replacement} of $M$.
\end{dfn}

The next construction recapitulates the proof of
\thmcite[A.6]{TNkPTh}; it is one of many ways to construct a
semi-flat-cotorsion replacement, \rmkref{HuDing} presents another.

\begin{con}
  \label{con:repl}
  Let $M$ be an $R$-complex and consider a semi-projective resolution
  $\qisdef{\pi^M}{P^M}{M}$. By \fctref{gil} there is an exact sequence
  \begin{equation*}
    0 \lra P^M \xra{\gre^M} C^M \lra A^M \lra 0
  \end{equation*}
  where $C^M$ is semi-cotorsion and $A^M$ is acyclic and semi-flat.
  It follows from \fctref{3-sf} that $C^M$ is semi-flat and isomorphic
  to $P^M$ and hence to $M$ in the derived category, so $C^M$ is a
  semi-flat-cotorsion replacement of $M$.
\end{con}

The Gorenstein flat-cotorsion dimension is defined based on
semi-flat-cotorsion replacements. Below we collect some technical
results for later use; the first one is about comparison of
semi-flat-cotorsion replacements.

\begin{prp}
  \label{prp:compare}
  Let $M$ be an $R$-complex. If $W$ and $W'$ are semi-flat-cotorsion
  replacements of $M$, then there is a homotopy equivalence
  $W \to W'$.
\end{prp}

\begin{prf*}
  Let $P \qra M$ be a semi-projective resolution; there are
  quasi-isomorphisms $W \qla P \qra W'$, see
  \cite[1.4.P]{LLAHBF91}. By \prpref{cothty} there is a
  quasi-isomorphism $W \qra W'$, and by \prpref{cot-hty} it is a
  homotopy equivalence.
\end{prf*}

\begin{rmk}
  \label{rmk:res}
  While, say, a semi-projective resolution is a quasi-isomorphism
  between complexes, there may not be a quasi-isomorphism between a
  complex and its semi-flat-cotorsion replacement. In certain cases,
  though, such maps do exist.

  Let $M$ be a semi-cotorsion complex. By \fctref{gil} there is an
  exact sequence of complexes $0 \lra C \lra F \xra{\varpi} M \lra 0$
  with $F$ semi-flat and $C$ semi-cotorsion and acyclic. It follows
  from \lemref{3-cot} that $F$ is semi-flat-cotorsion and, since $C$
  is acyclic the morphism $\varpi$ is a quasi-isomorphism.

  For a cotorsion module $M$, a semi-flat-cotorsion resolution
  $F \qra M$ can be constructed by taking successive flat covers,
  which have cotorsion kernels; in particular, one has $F_i=0$ for
  $i <0$.
\end{rmk}

The next result is a Schanuel's lemma for semi-flat-cotorsion
replacements.

\begin{lem}
  \label{lem:Schanuel}
  Let $M$ be an $R$-complex and $W$ and $W'$ be semi-flat-cotorsion
  replacements of $M$. For every $n \in \ZZ$ there exist
  flat-cotorsion $R$-modules $V$ and $V'$ with
  $\Co[n]{W} \oplus V \is \Co[n]{W'} \oplus V'$.
\end{lem}

\begin{prf*}
  By \prpref{compare} there is a homotopy equivalence
  $\mapdef{\gra}{W}{W'}$. The complex $\Cone{\gra}$ is contractible
  and semi-flat-cotorsion by \lemref{3-cot} and \fctref{3-sf}. It
  follows that every cycle module $\Cy[n]{\Cone{\gra}}$ is
  flat-cotorsion. The soft truncated morphism
  $\mapdef{\Tsa{n}{\gra}}{\Tsa{n}{W}}{\Tsa{n}{W'}}$ is also a homotopy
  equivalence, so $\Cone{(\Tsa{n}{\gra})}$, i.e.\ the complex
  \begin{equation*}
    0 \lra \Co[n]{W} \lra \Co[n]{W'} \oplus W_{n-1}
    \lra W'_{n-1} \oplus W_{n-2} \xra{\dif[n-1]{\Cone{\gra}}}
    W'_{n-2} \oplus W_{n-3} \lra \cdots
  \end{equation*}
  is contractible. Hence one has
  $\Co[n]{W} \oplus \Cy[n-1]{\Cone{\gra}} \is \Co[n]{W'} \oplus
  W_{n-1}$.
\end{prf*}

\begin{lem}
  \label{lem:trunc}
  Let $W$ be a semi-flat-cotorsion complex. For every $n\in \ZZ$ the
  truncated complex $\Thb{n}{W}$ is semi-flat-cotorsion.
\end{lem}

\begin{prf*}
  As $\Thb{n}{W}$ is a bounded below complex of flat modules it is
  semi-flat, and it is semi-cotorsion by \fctref{BCE}. Thus
  $\Thb{n}{W}$ is semi-flat-cotorsion.
\end{prf*}

\section{Gorenstein flat-cotorsion modules}

\noindent
Recall from \dfncite[4.3 and Prop.\ 1.3]{CETa} that a \emph{totally
  acyclic complex of flat-cotorsion modules} is an acyclic complex $T$
of flat-cotorsion modules such that the complexes $\Hom{T}{W}$ and
$\Hom{W}{T}$ are acyclic for every flat-cotorsion $R$-module $W$. An
$R$-module $G$ is called \emph{Gorenstein flat-cotorsion} if there
exists a totally acyclic complex $T$ of flat-cotorsion $R$-modules
with $\Cy[0]{T} = G$.

\begin{rmk}
  \label{rmk:cot-cy}
  To show that an acyclic complex $T$ of flat-cotorsion $R$-modules is
  totally acyclic it suffices to verify that $\Hom{T}{W}$ is acyclic
  for every flat-cotorsion $R$-module $W$; the acyclicity of
  $\Hom{W}{T}$ is automatic as the cycles in $T$ are cotorsion by
  \fctref{BCE}.
\end{rmk}

\begin{lem}
  \label{lem:enough}
  An $R$-module $G$ is Gorenstein flat-cotorsion if and only if the
  following conditions are satisfied.
  \begin{rqm}
  \item $G$ is cotorsion.
  \item $\Ext{i \ge 1}{G}{W}=0$ holds for every flat-cotorsion
    $R$-module $W$.
  \item There exists a complex $T = T_0 \to T_{-1} \to \cdots$ of
    flat-cotorsion $R$-modules and an injective quasi-isomorphism
    $\mapdef{\gre}{G}{T}$ such that $\Hom{\gre}{W}$ is a
    quasi-isomorphism for every flat-cotorsion $R$-module $W$.
  \end{rqm}
\end{lem}

\begin{prf*}
  The ``only if'' is in view of \rmkref{cot-cy} clear from the
  definition of Gorenstein flat-cotorsion modules. For the ``if''
  recall from \rmkref{res} that since $G$ is cotorsion, there is a
  semi-flat-cotorsion complex $T'$ and a surjective quasi-isomorphism
  $\qisdef{\pi}{T'}{G}$. Splicing $\Shift{T'}$ with the complex $T$
  yields an acyclic complex $T''$ of flat-cotorsion modules with
  $\Cy[0]{T''} = G$. For every flat-cotorsion module $W$ one has, by
  way of \fctref{rep}, that
  $$\H[i]{\Hom{T''}{W}} = \H[i+1]{\Hom{T'}{W}} = \Ext{-i-1}{G}{W} =0$$
  for $i \le -2$, and for $i\ge 0$ one has
  $\H[i]{\Hom{T''}{W}} = \H[i]{\Cone{\Hom{\gre}{W}}} = 0$. Finally
  $\H[-1]{\Hom{T''}{W}}=0$ holds as one has
  \begin{equation*}
    \Im{(\Hom{\gre_0}{W})} = \Hom{G}{W} = \Ker{(\Hom{\partial_1^{T'}}{W})}\:,
  \end{equation*}
  where the first equality holds as $\Hom{\gre}{W}$ is a
  quasi-isomorphism and the second holds by left exactness of
  $\Hom{-}{W}$.
\end{prf*}

\begin{prp}
  \label{prp:oplus}
  The class of Gorenstein flat-cotorsion $R$-modules is closed under
  finite direct sums and direct summands.
\end{prp}

\begin{prf*}
  To see that the class is closed under finite direct sums, let $T$
  and $T'$ be totally acyclic complexes of flat-cotorsion
  $R$-modules. The direct sum $T \oplus T'$ is a totally acyclic
  complex of flat-cotorsion modules with
  $\Cy[0]{T\oplus T'} = \Cy[0]{T} \oplus \Cy[0]{T'}$.

  Assume that $G$ and $G'$ are $R$-modules such that $G \oplus G'$ is
  Gorenstein flat-cotorsion. By additivity of Ext, it follows from
  \lemref{enough} that $G$ and $G'$ are cotorsion with
  $\Ext{i \ge 1}{G}{W}= 0 = \Ext{i \ge 1}{G'}{W}$ for every
  flat-cotorsion $R$-module $W$. It remains to verify part (3) of
  \lemref{enough} for, say, the module $G$.

  By definition there is an exact sequence of $R$-modules,
  \begin{equation*}
    0 \lra G\oplus G' \lra T_0 \lra G'' \lra 0\:,
  \end{equation*}
  with $T_0$ flat-cotorsion and $G''$ Gorenstein flat-cotorsion, in
  particular cotorsion. Consider the push-out diagram with exact rows
  and columns:
  \begin{equation*}
    \xymatrix@=1.5pc{
      & 0\ar[d] & 0\ar[d]\\
      & G\ar[d] \ar@{=}[r] & G \ar[d]\\
      0 \ar[r] & G\oplus G' \ar[d] \ar[r]
      & T_0 \ar[d]\ar[r] & G'' \ar@{=}[d] \ar[r] & 0 \\
      0 \ar[r] & G' \ar[d]\ar[r] & X_{-1} \ar[d] \ar[r] & G'' \ar[r] & 0 \\
      & 0 & 0 }
  \end{equation*}
  Let $W$ be a flat-cotorsion $R$-module.  As
  $\Ext{i \ge 1}{G'}{W} = 0 = \Ext{i \ge 1}{G''}{W}$, again by
  \lemref{enough}, one has $\Ext{i \ge 1}{X_{-1}}{W}=0$. Thus the
  exact sequence
  \begin{equation}
    \label{1}
    0 \lra G \lra T_0 \lra X_{-1} \lra 0
  \end{equation}
  is $\Hom{-}{W}$-exact for every flat-cotorsion $R$-module $W$.

  Next interchange the roles of $G'$ and $G$ and consider the diagram
  \begin{equation*}
    \xymatrix@=1.5pc{
      & 0\ar[d] & 0\ar[d]\\
      & G'\ar[d] \ar@{=}[r] & G' \ar[d]\\
      0 \ar[r] & G\oplus G' \ar[d] \ar[r] & T_0 \ar[d]\ar[r]
      & G'' \ar@{=}[d] \ar[r] & 0 \\
      0 \ar[r] & G \ar[d]\ar[r] & X'_{-1} \ar[d] \ar[r] & G'' \ar[r] & 0 \\
      & 0 & 0 }
  \end{equation*}
  with exact rows and columns. As above, $\Ext{i \ge 1}{X'_{-1}}{W}=0$
  holds for every flat-cotorsion $R$-module $W$. In particular, the
  exact sequence
  \begin{equation}
    \label{2}
    0 \lra G' \lra T_0 \lra X'_{-1} \lra 0
  \end{equation}
  is $\Hom{-}{W}$-exact for every flat-cotorsion $R$-module $W$.  The
  direct sum of the sequences (1) and (2) makes up the upper row in
  the commutative diagram
  \begin{equation*}
    \xymatrix@=1.5pc{
      0 \ar[r] & G\oplus G' \ar@{=}[d] \ar[r]
      & T_0\oplus T_0 \ar[d]_{\gra}\ar[r]
      & X_{-1}\oplus X'_{-1} \ar[d] \ar[r] & 0 \\
      0 \ar[r] & G\oplus G' \ar[r] & T_0 \ar[r] & G'' \ar[r] & 0 }
  \end{equation*}
  where $\gra\colon T_0\oplus T_0 \to T_0$ is the epimorphism given by
  $\gra(x,y)=x+y$.  By the Snake Lemma one gets the following
  commutative diagram with exact rows and columns
  \begin{equation*}
    \xymatrix@=1.5pc{
      & &0\ar[d] & 0\ar[d]\\
      & &T_0\ar[d] \ar@{=}[r] & T_0 \ar[d]\\
      0 \ar[r] & G\oplus G' \ar@{=}[d] \ar[r] & T_0\oplus T_0
      \ar[d]_{\gra}\ar[r] & X_{-1}\oplus X'_{-1} \ar[d] \ar[r] & 0 \\
      0 \ar[r] & G\oplus G' \ar[r] & T_0 \ar[d] \ar[r]
      & G'' \ar[d]\ar[r] & 0 \\
      & & 0 & 0  }
  \end{equation*}
  As $T_0$ is flat-cotorsion, the right-hand column splits, so there
  is an isomorphism $X_{-1}\oplus X'_{-1} \cong T_0 \oplus G''$. As
  the class of Gorenstein flat-cotorsion modules is closed under
  finite direct sums, it follows that $X_{-1}\oplus X'_{-1}$ is
  Gorenstein flat-cotorsion.

  Applying the same process to $X_{-1}\oplus X'_{-1}$ one gets exact
  sequences
  \begin{align*}
    &0 \lra X_{-1} \lra T_{-1} \lra X_{-2} \lra 0\quad\text{ and}\\
    &0 \lra X'_{-1} \lra T_{-1} \lra X'_{-2} \lra 0
  \end{align*}
  of $R$-modules where $T_{-1}$ is flat-cotorsion and
  $X_{-2}\oplus X'_{-2}$ is Gorenstein flat-cotorsion. Moreover, both
  sequences are $\Hom{-}{W}$-exact for every flat-cotorsion $R$-module
  $W$.  Continuing this process, one gets an exact sequence
  \begin{equation*}
    0 \lra G \lra T_0 \lra T_{-1} \lra T_{-2} \lra \cdots
  \end{equation*}
  with each $T_{i}$ flat-cotorsion, and the sequence is
  $\Hom{-}{W}$-exact for every flat-cotorsion $R$-module $W$.
\end{prf*}

\begin{prp}
  \label{prp:ses}
  Let $0 \to G' \to G \to G'' \to 0$ be an exact sequence of cotorsion
  $R$-modules.
  \begin{prt}
  \item If $G''$ is Gorenstein flat-cotorsion, then $G$ is Gorenstein
    flat-cotorsion if and only if $G'$ is Gorenstein flat-cotorsion.
  \item If $G'$ and $G$ are Gorenstein flat-cotorsion, then $G''$ is
    Gorenstein flat-cotorsion if and only if $\Ext{1}{G''}{W}=0$ holds
    for every flat-cotorsion $R$-module $W$.
  \end{prt}
\end{prp}

\begin{prf*}
  (a): Let $0 \lra G' \xra{\gra'} G \xra{\gra} G'' \lra 0$ be an exact
  sequence of cotorsion $R$-modules. If $G'$ and $G''$ are Gorenstein
  flat-cotorsion, then it follows from \prpcite[4.2 and
  Lem.~2.10]{CETa} that $G$ is Gorenstein flat-cotorsion.

  Assume now that $G$ and $G''$ are Gorenstein flat-cotorsion. By
  definition there exist exact sequences of $R$-modules
  \begin{equation*}
    0 \lra G \xra{\widetilde{\varepsilon}} \widetilde{T} \lra
    C \lra 0 \qqand
    0 \lra G'' \xra{\varepsilon''} T'' \lra \Coker{\varepsilon''} \lra 0
  \end{equation*}
  with $\widetilde{T}$ and $T''$ flat-cotorsion and $C$ and
  $\Coker{\varepsilon''}$ Gorenstein flat-cotorsion. One has
  $\Ext{1}{C}{T''}=0$ by \lemref{enough}, and so there is a
  homomorphism $\mapdef{\widetilde{\tau}}{\widetilde{T}}{T''}$ such
  that $\widetilde{\tau}\widetilde{\varepsilon} =\varepsilon''\alpha$
  holds. The module $T= \widetilde{T} \oplus T''$ is flat-cotorsion,
  and the surjective homomorphism
  $\mapdef{\tau=(\widetilde{\tau},1_{T''})}{T}{T''}$ is split. As the
  class of flat-cotorsion modules is closed under direct summands, the
  module $T' =\Ker{\tau}$ is flat-cotorsion. Let
  $\mapdef{\varepsilon}{G}{T}$ be the composite
  $ G \xra{\widetilde{\varepsilon}} \widetilde{T} \hookrightarrow T$;
  it is an injective homomorphism with $\Coker\varepsilon=C\oplus T''$
  Gorenstein flat-cotorsion, and
  $\tau\varepsilon = \widetilde{\tau}\widetilde{\varepsilon} =
  \varepsilon''\alpha$. Let $\tau'$ be the inclusion of $T'$ into $T$
  and $\mapdef{\varepsilon'}{G'}{T'}$ the injective homomorphism
  induced by $\varepsilon$. One now has a commutative diagram with
  exact rows
  \begin{equation*}
    \xymatrix@=1.5pc{
      0 \ar[r] & G' \ar[r]^-{\alpha'}\ar[d]_{\varepsilon'}
      & G \ar[r]^-{\alpha}\ar[d]_{\varepsilon} & G'' \ar[r]\ar[d]_{\varepsilon''}
      & 0\\
      0 \ar[r] & T' \ar[r]^-{\tau'} & T \ar[r]^-{\tau} & T'' \ar[r] & 0}
  \end{equation*}
  Consider the induced exact sequence
  $0 \to \Coker\varepsilon'\to \Coker\varepsilon \to
  \Coker\varepsilon''\to 0$. As $G'$ and $T'$ are cotorsion, so is the
  module $\Coker\varepsilon'$. Since $\Coker\varepsilon$ and
  $\Coker\varepsilon''$ are Gorenstein flat-cotorsion, one has
  $\Ext{1}{\Coker\varepsilon'}{W}=0$ for all flat-cotorsion
  $R$-modules $W$ by \lemref{enough}. It follows that the sequence
  \begin{equation*}
    0 \lra G'\xra{\varepsilon'} T' \lra \Coker\varepsilon' \lra 0
  \end{equation*}
  is $\Hom{-}{W}$-exact for all flat-cotorsion $R$-modules
  $W$. Repeating this process, one sees that the module $G'$ satisfies
  condition (3) in \lemref{enough}.  Moreover, $G'$ satisfies
  condition (1) by assumption, and as $G$ and $G''$ satisfy (2), so
  does $G'$.

  (b): The ``only if'' part follows from \lemref{enough}. For the
  converse, recall that there exists an exact sequence
  $0 \to G' \to T' \to G''' \to 0$ with $T'$ flat-cotorsion and $G'''$
  Gorenstein flat-cotorsion. Consider the push-out diagram
  \begin{equation}
    \label{eq:pushout-G}
    \begin{gathered}
      \xymatrix@=1.5pc{ {} & 0 \ar[d] & 0 \ar[d] & {} & {}
        \\
        0 \ar[r] & G' \ar[d] \ar[r] & G \ar[d] \ar[r] & G'' \ar@{=}[d]
        \ar[r] & 0
        \\
        0 \ar[r] & T' \ar[d] \ar[r] & X \ar[d] \ar[r] & G'' \ar[r] & 0
        \\
        {} & G''' \ar[d] \ar@{=}[r] & G''' \ar[d] & {} & {}
        \\
        {} & 0 & 0 & {} & {} }
    \end{gathered}
  \end{equation}
  with exact rows and columns. The second column and part (a) show
  that $X$ is Gorenstein flat-cotorsion. As $\Ext{1}{G''}{T'}=0$ holds
  by assumption, the second row splits, whence $G''$ is Gorenstein
  flat-cotorsion by \prpref{oplus}.
\end{prf*}

\section{Gorenstein flat-cotorsion dimension}
\noindent
We now turn to defining the advertised homological dimension, and
subsequently prove that it behaves as one would expect for all
associative rings.

\begin{dfn}
  \label{dfn:Gfcd}
  Let $M$ be an $R$-complex. The \emph{Gorenstein flat-cotorsion
    dimension} of $M$, written $\Gfcd{M}$, is defined as
  \begin{equation*}
    \Gfcd{M} \deq \inf%
    \left\{n\in\ZZ
      \:\left|\:
        \begin{gathered}
          \text{There is a semi-flat-cotorsion replacement}\\[-3pt]
          \text{$W$ of $M$ with $\H[i]{W}=0$ for all $i > n$ and} \\[-3pt]
          \text{ with $\Co[n]{W}$ Gorenstein flat-cotorsion}
        \end{gathered}
      \right.
    \right\}
  \end{equation*}
  with $\inf\varnothing = \infty$ by convention.  We say that
  $\Gfcd{M}$ is \emph{finite} if $\Gfcd{M}<\infty$.
\end{dfn}

\begin{ipg}
  Let $M$ be an $R$-complex. For every semi-flat-cotorsion replacement
  $W$ of $M$ one has $\H{W} \is \H{M}$; the next (in)equalities are
  hence immediate from the definition,
  \begin{equation*}
    \Gfcd{M} \ge \sup{M} \qqand \Gfcd{\Shift[s]{M}} = \Gfcd{M} + s
    \ \text{ for every $s\in\ZZ$\;.}
  \end{equation*}
  Moreover, one has $\Gfcd{M} = -\infty$ if and only if $M$ is
  acyclic.
\end{ipg}

\begin{lem}
  \label{lem:any-Co}
  Let $M$ be an $R$-complex. For every semi-flat-cotorsion replacement
  $W$ of $M$ and every $n \ge \Gfcd{M}$ the module $\Co[n]{W}$ is
  Gorenstein flat-cotorsion.
\end{lem}

\begin{prf*}
  Assume that $\Gfcd{M} = g$ holds for some integer $g$. By assumption
  there exists a semi-flat-cotorsion replacement $W$ of $M$ with
  $\Co[g]{W}$ Gorenstein flat-cotorsion and $g \ge \sup{W}$. By
  induction it follows from \prpref{ses}(a) that $\Co[n]{W}$ is
  Gorenstein flat-cotorsion for every $n \ge g$. Let $W'$ be any
  semi-flat-cotorsion replacement of $M$. It follows from
  \lemref{Schanuel} and \prpref{oplus} that $\Co[n]{W'}$ is Gorenstein
  flat-cotorsion for every $n \ge g$.
\end{prf*}

\begin{prp}
  \label{prp:rep}
  Let $M$ be an $R$-complex.
  \begin{prt}
  \item For every bounded $R$-complex $C$ of flat-cotorsion modules
    one has
    \begin{equation*}
      \inf{\RHom{M}{C}} \dge \inf\setof{i\in\ZZ}{C_i\ne 0} - \Gfcd{M}\:.
  \end{equation*}
  \item
    For every complex $F$ of finite flat dimension one has
    \begin{equation*}
      \sup{\RHom{F}{M}} \dle \Gfcd{M} - \inf{F}\:.
    \end{equation*}
  \end{prt}
\end{prp}

\begin{prf*}
  Both assertions are trivial if $M$ is acyclic so assume that it is
  not. The assertions are also trivial if $\Gfcd{M}$ is infinite, so
  assume that it is finite and set $n = \Gfcd{M}$. Let $W$ be a
  semi-flat-cotorsion replacement of $M$. There is a short exact
  sequence $0 \to W' \to W \to \Tsa{n}{W} \to 0$; by \lemref{any-Co}
  the module $\Co[n]{W}$ is Gorenstein flat-cotorsion, so it follows
  from \prpref{ses}(a) that $W'$ is a complex of Gorenstein
  flat-cotorsion modules, concentrated in degrees $n$ and
  above. Further, $W'$ is acyclic, as the map $W \to \Tsa{n}{W}$ is a
  quasi-isomorphism. By \lemref{any-Co} the cycle modules in $W'$ are
  Gorenstein flat-cotorsion, as one has $\Cy[i]{W'} \is \Co[i+1]{W}$
  for all $i\ge n$.

  (a): For every flat-cotorsion module $C'$ the complex $\Hom{W'}{C'}$
  is acyclic; see \lemref{enough}. As $C$ is bounded above, $\Hom{W'}{C}$ is acyclic by \lemcite[2.5]{CFH-06}. Thus
  \fctref{rep} yields an isomorphism
  $\RHom{M}{C} \qis \Hom{\Tsa{n}{W}}{C}$ in the derived category. A
  standard argument---see for example the proof of
  \thmcite[3.1]{CFH-06}---yields $\Hom{\Tsa{n}{W}}{C}_i=0$ for
  $i < \inf\setof{i\in\ZZ}{C_i \ne 0} - n$.

  (b) We can assume that the homology of $F$ is bounded below,
  otherwise the inequality is trivial. By \rmkcite[1.7 and
  Thm.~2.4.F]{LLAHBF91} we can further assume that $F$ is a bounded
  complex of flat modules with $F_i = 0$ for $i < \inf{F}$.  For every
  $i \in \ZZ$, the complex $\Hom{F_i}{W'}$ is acyclic, as the cycle
  modules in $W'$ are cotorsion. As $F$ is bounded, it follows from
  \lemcite[2.4]{CFH-06} that $\Hom{F}{W'}$ is acyclic. By \fctref{rep}
  there is thus an isomorphism $\RHom{F}{M} \qis \Hom{F}{\Tsa{n}{W}}$
  in the derived category.  A standard argument now yields
  $\Hom{F}{\Tsa{n}{W}}_i=0$ for $i > -\inf{F} + n$.
\end{prf*}

\begin{thm}
  \label{thm:Gfcd}
  Let $M$ be an $R$-complex and $n$ be an integer. If\, $\Gfcd{M}$ is
  finite and $n \ge \sup{M}$, then the following conditions are
  equivalent.
  \begin{eqc}
  \item $\Gfcd{M} \le n$.
  \item $-\inf{\RHom{M}{C}} \le n$ for every cotorsion $R$-module $C$
    of finite flat dimension.
  \item $-\inf{\RHom{M}{C}} \le n$ for every flat-cotorsion $R$-module
    $C$.
  \item $\Ext{n+1}{M}{C} =0$ for every cotorsion $R$-module $C$ of
    finite flat dimension.
  \item For every semi-flat-cotorsion replacement $W$ of $M$ the
    module $\Co[m]{W}$ is Gorenstein flat-cotorsion for every
    $m \ge n$.
  \end{eqc}
\end{thm}

The proof relies on a technical lemma which follows at the end of the
section.

\begin{prf*}
  The statement is trivial for an acyclic complex, so assume that $M$
  is not acyclic. The implications \proofofimp[]{ii}{iii},
  \proofofimp[]{ii}{iv}, and \proofofimp[]{v}{i} are clear. The
  implication \proofofimp[]{i}{v} follows from \lemref{any-Co}.

  \proofofimp{i}{ii} Apply \prpref{rep}(a) to a finite flat resolution
  of $C$ by flat-cotorsion modules; see \rmkref{res}.

  \proofofimp{iii}{i} Let $W$ be a semi-flat-cotorsion replacement of
  $M$. As $g=\Gfcd{M}$ is finite, it follows from \lemref{any-Co} that
  the module $\Co[m]{W}$ is Gorenstein flat-cotorsion for every
  integer $m \ge g$. Thus it is enough to show that $n \ge
  g$ holds. Assume towards a contradiction that $n<g$ holds. It now follows
  from the assumption $n \ge \sup{M}$ that the module $\Co[g-1]{W}$ is
  not Gorenstein flat-cotorsion.  By \lemref{trunc} the complex
  $\Shift[1-g]{\Thb{g-1}{W}}$ is a semi-flat-cotorsion replacement of
  the module $\Co[g-1]{W}$. Let $C$ be a flat-cotorsion $R$-module; in
  view of \fctref{rep} one has
  \begin{align*}
    \Ext{1}{\Co[g-1]{W}}{C} 
    & \is \H[-1]{\RHom{\Co[g-1]{W}}{C}} \\
    & = \H[-1]{\Hom{\Shift[1-g]{\Thb{g-1}{W}}}{C}} \\
    & = \H[-g]{\Hom{\Thb{g-1}{W}}{C}} \\
    & = \H[-g]{\RHom{M}{C}} \\
    & = 0\:.
  \end{align*}
  It now follows by an application of \prpref{ses}(b) to the short
  exact sequence $0 \to \Co[g]{W} \to W_{g-1} \to \Co[g-1]{W} \to 0$
  that $\Co[g-1]{W}$ is Gorenstein flat-cotorsion; a contradiction.

  \proofofimp{iv}{i} Let $W$ be a semi-flat-cotorsion replacement of
  $M$. As in the argument for \proofofimp[]{iii}{i}, set $g=\Gfcd{M}$,
   assume towards a contradiction that $n<g$ holds, and notice that
  the module $\Co[n]{W}$ is not Gorenstein flat-cotorsion.  There is
  an exact sequence
  $0 \to \Co[g]{W} \to W_{g-1} \to \cdots \to W_n \to \Co[n]{W} \to
  0$, so by \lemref{hh} there is an exact sequence of cotorsion
  $R$-modules
  \begin{equation}
    \tag{1}
    0 \lra K \lra G \lra \Co[n]{W} \lra 0
  \end{equation}
  where $G$ is Gorenstein flat-cotorsion.
  The complex $\Shift[-n]{\Thb{n}{W}}$ is a semi-flat-cotorsion
  replacement of the module $\Co[n]{W}$, see \lemref{trunc}, so as
  above one has
  \begin{align*}
    \Ext{1}{\Co[n]{W}}{K} 
    & \is \H[-1]{\RHom{\Co[n]{W}}{K}} \\
    & = \H[-1]{\Hom{\Shift[-n]{\Thb{n}{W}}}{K}} \\
    & = \H[-(n+1)]{\Hom{\Thb{n}{W}}{K}} \\
    & \is \Ext{n+1}{M}{K} \\
    & = 0      \:.
  \end{align*}
  It follows that the sequence (1) splits, which by \prpref{oplus}
  implies that $\Co[n]{W}$ is Gorenstein flat-cotorsion; a
  contradiction.
\end{prf*}

\begin{rmk}
  \label{rmk:cotorsionmodules}
  As one can surmise from the statement of the next lemma, not every
  $R$-module has a flat-cotorsion resolution---an example of such a
  module is provided in \exacite[3.11]{TNkPTh}---but per \rmkref{res}
  cotorsion modules do. Thus, for a cotorsion $R$-module $C$ one has
  $\Gfcd{C} \le n$ if and only if there is an exact sequence
  $0 \to G \to W_{n-1} \to \cdots \to W_0 \to C \to 0$ where the
  modules $W_i$ are flat-cotorsion and $G$ is Gorenstein
  flat-cotorsion. In particular, a cotorsion $R$-module $C$ is
  Gorenstein flat-cotorsion if and only if $\Gfcd{C} = 0$ holds.
\end{rmk}

The proof of the next lemma is adapted from Holm \cite{HHl04a}.

\begin{lem}
  \label{lem:hh}
  Let $M$ be an $R$-module. If there exists an exact sequence
  \begin{equation*}
    0 \lra G \lra W_{g-1} \lra \cdots \lra W_0 \lra M \lra 0
  \end{equation*}
  in which each module $W_i$ is flat-cotorsion and $G$ is Gorenstein
  flat-cotorsion, then there is an exact sequence of cotorsion
  $R$-modules $0 \to K \to G' \to M \to 0$ and a quasi-isomorphism
  $V \qra K$ where $G'$ is Gorenstein flat-cotorsion and $V$ is a
  semi-flat-cotorsion complex concentrated in degrees $g-1,\ldots,0$.
\end{lem}

\begin{prf*}
  As $G$ is Gorenstein flat-cotorsion, the defining totally acyclic
  complex yields an exact sequence
  \begin{equation*}
    0 \lra G \lra W'_{g-1} \lra \cdots \lra W'_0 \lra H \lra 0
  \end{equation*}
  where each module $W_i'$ is flat-cotorsion and $H$ is Gorenstein
  flat-cotorsion. Moreover, this sequence is $\Hom{-}{W''}$-exact for
  every flat-cotorsion $R$-module $W''$. It follows that the identity
  on $G$ lifts to a morphism of complexes
  \begin{equation*}
    \xymatrix@=1.5pc{
      0 \ar[r] &
      W'_{g-1} \ar[r]\ar[d] & \cdots \ar[r] & W'_0 \ar[d] \ar[r]
      & H \ar[d] \ar[r] & 0\\
      0 \ar[r] & W_{g-1} \ar[r] & \cdots \ar[r] & W_0 \ar[r] & M \ar[r] & 0
    }
  \end{equation*}
  It is a quasi-isomorphism, so the mapping cone
  \begin{equation*}
    0 \lra W'_{g-1} \lra W'_{g-2}\oplus W_{g-1} \lra \cdots
    \lra  W'_{0}\oplus W_{1}
    \lra H \oplus W_0 \lra M \lra 0
  \end{equation*}
  is acyclic. The module $G' = H \oplus W_0$ is Gorenstein
  flat-cotorsion by \prpref{oplus}. The truncated complex
  $V = 0 \to W'_{g-1} \to \cdots \to W'_{0}\oplus W_{1} \to 0$ is
  semi-flat-cotorsion, and with $K = \Ker{(G' \to M)}$ the canonical
  map $V \to K$ is a quasi-isomorphism. The modules $K$ and $M$ are
  both cotorsion as cotorsion modules are closed under cokernels of
  monomorphisms.
\end{prf*}

\begin{rmk}
  As observed in \prpcite[4.2]{CETa}, a Gorenstein flat-cotorsion
  module is the same as a right $\mathsf{Flat}(R)$-Gorenstein module
  in the sense of \dfncite[2.1]{CETa}. Given a sufficiently
  well-behaved cotorsion pair $(\mathsf{U},\mathsf{V})$ in the
  category $\mathsf{Mod}(R)$ of $R$-modules, one can adapt the
  definitions and arguments in Sections 1--4 of this paper to develop
  a right \textsf{U}-Gorenstein dimension theory for $R$-complexes,
  and one can dualize them to develop a left \textsf{V}-Gorenstein
  dimension theory.

  Applied to the cotorsion pair
  $(\mathsf{Prj}(R),\mathsf{Mod}(R))$ the right dimension
  is the classic Gorenstein projective dimension and the left
  dimension is trivial; it detects the homological infimum of a
  complex; see \cite[Exmpl.~2.5]{CETa}. Dually, the cotorsion pair
  $(\mathsf{Mod}(R),\mathsf{Inj}(R))$ yields a
  right dimension that detects the homological supremum of a complex
  while the left dimension is the classic Gorenstein injective
  dimension.

  For the cotorsion pair $(\mathsf{Flat}(R),\mathsf{Cot}(R))$ the
  right $\mathsf{Flat}(R)$-Gorenstein dimension is the one developed
  in this paper under a different name. A left
  $\mathsf{Cot}(R)$-Gorenstein module is simply a flat-cotorsion
  module, see again \prpcite[4.2]{CETa}, so the left
  $\mathsf{Cot}(R)$-Gorenstein dimension of a complex $M$ is in this
  case the least $n$ such that there is a semi-flat-cotorsion
  replacement $F$ of $M$ with $F_i=0$ for all $i < -n$.
\end{rmk}

\section{Comparison to the Gorenstein flat dimension}
\label{sec:compare}

\noindent
In this section, we put the newly minted Gorenstein flat-cotorsion
dimension into context, showing that in most settings it is nothing
but an avatar of the classic Gorenstein flat dimension. Recall that a
complex of flat $R$-modules is called \emph{\Ftac} if it is acyclic
and $\tp{E}{F}$ is acyclic for every injective $\Rop$-module $E$. An
$R$-module $G$ is called \emph{Gorenstein flat} if there exists an
\Ftac\ complex $F$ with $\Cy[0]{F} = G$.

A main result of \cite{JSrJSt20} is that the class of Gorenstein flat
modules is the left class in a hereditary cotorsion pair. From
\corcite[4.12]{JSrJSt20} one can, in particular, extract the following
statement.

\begin{fct}
  \label{fct:ss}
  The class of Gorenstein flat $R$-modules is closed under direct
  summands, extensions, and kernels of epimorphisms. Also, for every
  Gorenstein flat $R$-module $G$ one has $\Ext{i}{G}{W}=0$ for all
  $i\ge 1$ and every flat-cotorsion $R$-module $W$.
\end{fct}

The first step towards comparison of the dimensions is to compare
Gorenstein flat modules to Gorenstein flat-cotorsion modules.
\begin{thm}
  \label{thm:Gfc}
  Let $G$ be an $R$-module. If $G$ is Gorenstein flat and cotorsion,
  then it is Gorenstein flat-cotorsion. The converse holds if $R$ is
  right coherent.
\end{thm}

\begin{prf*}
  Assume that $G$ is Gorenstein flat and cotorsion; it satisfies
  condition (1) in \lemref{enough} and by \fctref{ss} also condition
  (2). Thus, it suffices to show that $G$ satisfies condition (3) in
  \lemref{enough}. By definition there is a short exact sequence
  \mbox{$0 \to G \to F \to G' \to 0$} of $R$-modules with $F$ flat and
  $G'$ Gorenstein flat. There is also an exact sequence
  $0 \to F \to W \to F' \to 0$ of $R$-modules with $W$ flat-cotorsion
  and $F'$ flat. Consider the push-out diagram
  \begin{equation*}
    \xymatrix@=1.5pc{
      & & 0 \ar[d] & 0 \ar[d] \\
      0 \ar[r] & G \ar@{=}[d] \ar[r] & F \ar[d] \ar[r] & G' \ar[d] \ar[r] & 0 \\
      0 \ar[r] & G \ar[r]^-{\gre_0} & W \ar[d] \ar[r] & X \ar[d] \ar[r] & 0 \\
      & & F'\ar[d] \ar@{=}[r] & F' \ar[d]\\
      & & 0 & 0}
  \end{equation*}
  with exact rows and columns. By \fctref{ss}, exactness of the
  right-hand column implies that $X$ is Gorenstein flat, and it
  follows that the middle row is $\Hom{-}{V}$-exact for every
  flat-cotorsion $R$-module $V$.  As $G$ and $W$ are cotorsion, so is
  $X$. Continuing this process, one gets an exact sequence
  \begin{equation}
    \label{eq:1}
    0 \lra G \xra{\gre_0} W_{-1} \lra W_{-2} \lra \cdots
  \end{equation}
  that is $\Hom{-}{V}$-exact for every flat-cotorsion $R$-module
  $V$. The modules $W_n$ are flat-cotorsion, so
  $W_{-1} \to W_{-2} \to \cdots$ is the desired complex, and $\gre_0$
  is the non-zero component of the desired quasi-isomorphism.

  If $R$ is right coherent, then the converse holds by
  \thmcite[5.2]{CETa}.
\end{prf*}

The definitions of Gorenstein flat dimension found in the literature
all agree with the following definition; see also the discussion in
\cite[Rmk.~5.13]{CKL-17}.

\begin{dfn}
  \label{dfn:Gfd}
  Let $M$ be an $R$-complex. The \emph{Gorenstein flat dimension} of
  $M$, written $\Gfd{M}$, is defined as
  \begin{equation*}
    \Gfd{M} \deq \inf%
    \left\{n\in\ZZ
      \:\left|\:
        \begin{gathered}
          \text{There exists a semi-flat replacement $F$ of $M$ with}\\[-3pt]
          \text{ $\H[i]{F}=0$ for all $i > n$ and $\Co[n]{F}$
            Gorenstein flat}
        \end{gathered}
      \right.
    \right\}
  \end{equation*}
  with $\inf\varnothing = \infty$ by convention.  We say that
  $\Gfd{M}$ is \emph{finite} if $\Gfd{M}<\infty$.
\end{dfn}

\noindent A crucial consequence of \fctref{ss} is \lemref{abc}
below. We provide a direct proof but note that it can also be deduced
by combining  results from \cite{LLn} and \cite{JSrJSt20}. One
can also obtain \lemref{abc} from arguments by Sather-Wagstaff,
Sharif, and White \cite{SSW-11}, still in view of \cite{JSrJSt20}.

\begin{lem}
  \label{lem:abc}
  Let $M$ be an $R$-module. If $M$ has finite Gorenstein flat
  dimension, then the following conditions are equivalent
  \begin{eqc}
  \item $M$ is Gorenstein flat.
  \item $\Ext{i}{M}{C} =0$ holds for all $i\ge 1$ and every
    flat-cotorsion $R$-module $C$.
  \item $\Ext{i}{M}{C} =0$ holds for all $i\ge 1$ and every
    cotorsion $R$-module $C$ of finite flat dimension.
  \item $\Ext{1}{M}{C} =0$ holds for  every
    cotorsion $R$-module $C$ of finite flat dimension.
\end{eqc}
\end{lem}

\begin{prf*}
  The implication \proofofimp[]{iii}{iv} is trivial and
  \proofofimp[]{i}{ii} follows from \fctref{ss}.

  \proofofimp{ii}{iii} A cotorsion module $C$ of finite flat dimension
  has a bounded flat resolution by flat-cotorsion modules; see
  \rmkref{res}. The claim now follows by dimension shifting.

  \proofofimp{iv}{i} We first reduce to the case where $M$ is
  cotorsion. To this end, consider an exact sequence
  $0 \to M \to C^M \to F^M \to 0$ of $R$-modules, where $C^M$ is
  cotorsion and $F^M$ is flat.  By \fctref{ss} it suffices to prove
  that $C^M$ is Gorenstein flat. From the Horseshoe Lemma for
  projective resolutions and \fctref{ss} it follows that $C^M$ has
  finite Gorenstein flat dimension. We may now assume that $M$ is
  cotorsion.

  As $M$ has finite Gorenstein flat dimension, there is an exact
  sequence
  \begin{equation*}
    0 \lra G \lra W_{g-1} \lra \cdots \lra W_0 \lra M \lra 0
  \end{equation*}
  where the modules $W_i$ are flat-cotorsion and $G$ is Gorenstein
  flat. From the \Ftac\ complex defining $G$ one gets an exact
  sequence
  \begin{equation*}
    0 \lra G \lra F_{g-1} \lra \cdots \lra F_0 \lra H \lra 0
  \end{equation*}
  where each module $F_i$ is flat and $H$ is Gorenstein flat. By
  \fctref{ss} this sequence is $\Hom{-}{V}$-exact for every
  flat-cotorsion $R$-module $V$. As in the proof of \lemref{hh} one
  now gets an acyclic complex
  \begin{equation*}
    0 \lra F_{g-1} \lra F_{g-2}\oplus W_{g-1} \lra \cdots \lra  F_{0}\oplus W_{1}
    \lra H \oplus W_0 \lra M \lra 0\;.
  \end{equation*}
  The module $G' = H \oplus W_0$ is Gorenstein flat, and the module
  $K = \Ker{(G' \to M)}$ has finite flat dimension. There is a short
  exact sequence $0 \to K \to C^K \to F^K \to 0$ of $R$-modules with
  $C^K$ cotorsion and $F^K$ flat. It yields a pushout diagram
  \begin{equation}
    \label{eq:pushout-G2}
    \begin{gathered}
      \xymatrix@=1.5pc{ {} & 0 \ar[d] & 0 \ar[d] & {} & {}
        \\
        0 \ar[r] & K \ar[d] \ar[r] & G' \ar[d] \ar[r] & M \ar@{=}[d]
        \ar[r] & 0
        \\
        0 \ar[r] & C^K \ar[d] \ar[r] & X \ar[d] \ar[r] & M \ar[r] & 0
        \\
        {} & F^K \ar[d] \ar@{=}[r] & F^K \ar[d] & {} & {}
        \\
        {} & 0 & 0 & {} & {} }
    \end{gathered}
  \end{equation}
  with exact rows and columns. By \fctref{ss} the module $X$ is
  Gorenstein flat. The module $C^K$ is cotorsion of finite flat
  dimension, so by assumption the second row splits, which means that
  $M$ is Gorenstein flat in view of \fctref{ss}.
\end{prf*}

\begin{prp}
  \label{prp:abc}
  Let $M$ be an $R$-complex of finite Gorenstein flat dimension and
  \mbox{$n \ge \sup{M}$} an integer. One has $\Gfd{M} \le n$ if and
  only if $-\inf{\RHom{M}{C}} \le n$ holds for every cotorsion
  $R$-module $C$ of finite flat dimension.
\end{prp}

\begin{prf*}
  Assume that $\Gfd{M} \le n$ holds and let $F$ be a semi-flat
  replacement of $M$ as in \dfnref{Gfd}. As in the proof of
  \thmref{Gfcd} there is a short exact sequence
  $0 \to F' \to F \to \Tsa{n}{F} \to 0$. The complex $F'$ is acyclic
  and concentrated in degrees $n$ and above; by \fctref{ss} it is a
  complex of Gorenstein flat modules. For the same reason, the cycle
  modules in $F'$ are Gorenstein flat. Let $C$ be a cotorsion module
  of finite flat dimension. It follows from \lemref{abc} that
  $\Hom{F'}{C}$ is acyclic, so there is an isomorphism
  $\RHom{M}{C} \qis \Hom{\Tsa{n}{F}}{C}$ in the derived category.
  Evidently, one has $\Hom{\Tsa{n}{F}}{C}_i=0$ for $i < - n$.

  For the converse, set $g=\Gfd{M}$. One has to show that $n \ge g$
  holds.  Assume towards a contradiction that $n<g$ holds. It now
  follows from the assumption $n \ge \sup{M}$ that the module
  $\Co[g-1]{F}$ is not Gorenstein flat.  The complex
  $\Shift[1-g]{\Thb{g-1}{F}}$ is a semi-flat replacement of the module
  $\Co[g-1]{F}$ which, therefore, is a module of finite Gorenstein
  flat dimension $1$. Let $C$ be a cotorsion $R$-module of finite flat
  dimension; in view of \fctref{rep} one has
  \begin{align*}
    \Ext{1}{\Co[g-1]{F}}{C} 
    & \is \H[-1]{\RHom{\Co[g-1]{F}}{C}} \\
    & = \H[-1]{\Hom{\Shift[1-g]{\Thb{g-1}{F}}}{C}} \\
    & = \H[-g]{\Hom{\Thb{g-1}{F}}{C}} \\
    & = \H[-g]{\Hom{F}{C}} \\
    & = \H[-g]{\RHom{M}{C}} \\
    & = 0\:.
  \end{align*}
  It now follows from \lemref{abc} that $\Co[g-1]{F}$ is Gorenstein
  flat, a contradiction.
\end{prf*}

\begin{lem}
  \label{lem:Cy-cot}
  Let $C$ be a complex of cotorsion $R$-modules. For every
  $n \ge \sup{C}$ the module $\Co[n]{C}$ is cotorsion.
\end{lem}

\begin{prf*}
  Let $n \ge \sup{C}$. Splicing a shifted injective resolution of
  $\Co[n]{C}$ with the acyclic complex
  $\cdots \to C_{n+1} \to C_{n} \to \Co[n]{C} \to 0$ one gets an
  acyclic complex $X$ of cotorsion modules, and it follows from
  \fctref{BCE} that the cycle module $\Cy[n-1]{X} \is \Co[n]{C}$ is
  cotorsion.
\end{prf*}

The next theorem justifies the title of the paper.

\begin{thm}
  \label{thm:Gfd-Gfcd}
  Let $M$ be an $R$-complex. There is an inequality
  \begin{equation*}
    \Gfcd{M} \dle \Gfd{M}\:,
  \end{equation*}
  and equality holds if $M$ has finite Gorenstein flat dimension.
\end{thm}

\begin{prf*}
  It is enough to prove that every complex of finite Gorenstein flat
  dimension has finite Gorenstein flat-cotorsion dimension, then
  \thmref{Gfcd} and \prpref{abc} show that they agree.  Assume that
  $\Gfd{M} = n$ holds for some integer $n$. That is, there exists a
  semi-flat replacement $F$ of $M$ with $\Co[n]{F}$ Gorenstein
  flat. Consider an approximation $0 \to F \to C \to A \to 0$, where
  $A$ is acyclic and semi-flat and $C$ is semi-cotorsion; see
  \fctref{gil}. It follows that $C$ is a semi-flat-cotorsion
  replacement of $M$; see \fctref{3-sf}. As $A$ is acyclic there is an
  exact sequence $0 \to \Co[n]{F} \to \Co[n]{C} \to \Co[n]{A} \to
  0$. The module $\Co[n]{A}$ is flat and the class of Gorenstein flat
  modules is closed under extensions, see \fctref{ss}, so $\Co[n]{C}$
  is Gorenstein flat. By \lemref{Cy-cot} it is also cotorsion, so by
  \thmref{Gfc} it is Gorenstein flat-cotorsion.
\end{prf*}

\begin{cor}
  \label{cor:coh}
  Let $R$ be right coherent. For every $R$-complex $M$ one has
  \begin{equation*}
    \Gfd{M} = \Gfcd{M}\:.
  \end{equation*}
\end{cor}

\begin{prf*}
  It suffices in view of \thmref{Gfd-Gfcd} to prove that every complex
  of finite Gorenstein flat-cotorsion dimension has finite Gorenstein
  flat dimension, and that is immediate from the definitions and
  \thmref{Gfc}.
\end{prf*}

\begin{rmk}[a caveat]
  Every Gorenstein flat $R$-module $G$ has Gorenstein flat dimension
  $0$ and hence $\Gfcd{G}=0$ as well. Thus, a module of Gorenstein
  flat-cotorsion dimension $0$ need not be Gorenstein flat-cotorsion;
  any flat module that is not cotorsion exemplifies this; see also
  \rmkref{cotorsionmodules}.

  So how far are Gorenstein flat modules from being Gorenstein
  flat-cotorsion? \fctref{ss} can be harnessed to provide an answer:
  Let $G$ be a Gorenstein flat $R$-module. There are exact sequences
  of $R$-modules,
  \begin{equation*}
    0 \lra G \lra H \lra F \lra 0 \qqand 0 \lra H' \lra F' \lra G \lra 0 \;,
  \end{equation*}
  where $H$ and $H'$ are Gorenstein flat-cotorsion and $F$ and $F'$
  are flat. Indeed, the left-hand sequence exists with $H$ cotorsion
  and $F$ flat, and by \fctref{ss} the module $H$ is Gorenstein flat
  and, hence, Gorenstein flat-cotorsion by \thmref{Gfc}. Similarly,
  the right-hand sequence exists with $F'$ flat and $H'$ cotorsion,
  and by \fctref{ss} the module $H'$ is Gorenstein flat and, hence,
  Gorenstein flat-cotorsion by \thmref{Gfc}.
\end{rmk}

Recall that a noetherian (i.e.\ noetherian on both sides) ring is
called Iwanaga--Gorenstein if both injective dimensions $\id{R}$ and
$\id[\Rop]{R}$ are finite.

\begin{cor}
  Let $R$ be noetherian. The following conditions are equivalent.
  \begin{eqc}
  \item $R$ is Iwanaga--Gorenstein.
  \item Every $R$- and every $\Rop$-module has finite Gorenstein
    flat-cotorsion dimension.
  \item Every $R$- and every $\Rop$-complex with bounded above
    homology has finite Gorenstein flat-cotorsion dimension.
  \end{eqc}
\end{cor}

\begin{prf*}
  As $R$ is noetherian, it follows from \corref{coh} that the
  Gorenstein flat dimension coincides with the Gorenstein
  flat-cotorsion dimension over both $R$ and $\Rop$.  The assertion
  now follows from results of Enochs and Jenda \thmcite[12.3.1]{rha}
  and Iacob \thmcite[3.2]{AIc09}.
\end{prf*}

\begin{lem}
  \label{lem:fdlemma}
  Let $M$ be an $R$-complex and $n$ an integer. One has $\fd{M} \le n$
  if and only if $M$ has a semi-flat-cotorsion replacement $W$ with
  $W_i=0$ for all $i > n$.
\end{lem}

\begin{prf*}
  The ``if'' statement is evident. For the ``only if,'' let $W$ be a
  semi-flat-cotorsion replacement of $M$. By assumption the module
  $\Co[n]{W}$ is flat, see \thmcite[2.4.F]{LLAHBF91}, so there is an
  exact sequence of complexes of flat $R$-modules
  \begin{equation*}
    0 \lra W' \lra W \lra \Tsa{n}{W} \lra 0\:.
  \end{equation*}
  The sequence is degreewise pure exact, so it follows from
  \prpcite[6.2]{LWCHHl15} that the complex $\Tsa{n}{W}$ is
  semi-flat. As one has $n \ge \sup{M}$, the module $\Co[n]{W}$ is
  cotorsion by \lemref{Cy-cot}. It follows that
  $\Tsa{n}{W}$ consists of cotorsion $R$-modules, whence it is
  semi-cotorsion by \fctref{BCE}. Thus $\Tsa{n}{W}$ is a
  semi-flat-cotorsion replacement of~$M$.
\end{prf*}

\begin{thm}
  \label{thm:fd-Gfcd}
  Let $M$ be an $R$-complex. There is an inequality,
  \begin{equation*}
    \Gfcd{M} \dle \fd{M}\:,
  \end{equation*}
  and equality holds if $M$ has finite flat dimension.
\end{thm}

\begin{prf*}
  Without loss of generality assume that $n := \fd{M}$ is an
  integer. By \lemref{fdlemma} there is a semi-flat-cotorsion
  replacement $W$ of $M$ concentrated in degrees $n$ and below. In
  particular, one has $\Gfcd{M} \le n$. Set $m = \Gfcd{M}$ and assume
  towards a contradiction that $n > m$ holds. By \thmref{Gfcd} the
  module $\Co[m]{W}$ is Gorenstein flat-cotorsion, and as
  $m \ge \sup{M}$ holds there is an exact sequence of $R$-modules,
  \begin{equation*}
    0 \lra W_n \lra W_{n-1} \lra \cdots \lra W_m \lra \Co[m]{W} \lra 0\:.
  \end{equation*}
  It follows that $\Co[m]{W}$ has finite flat dimension; in
  particular, it has finite Gorenstein flat dimension. By
  \thmref{Gfd-Gfcd} one now has $\Gfd{\Co[m]{W}}=0$, which means that
  $\Co[m]{W}$ is a Gorenstein flat module. Now it follows from
  \fctref{ss} that the module $\Co[i]{W}$ is Gorenstein flat for
  $n > i \ge m$. It follows from \lemref{abc} that the following
  sequence splits
  \begin{equation*}
    0 \lra W_n \lra W_{n-1} \lra \Co[n-1]{W} \lra 0 \:,
  \end{equation*}
  whence $\Co[n-1]{W}$ is flat-cotorsion. So one has $\fd{M} \le n-1$,
  a contradiction.
\end{prf*}

\section{An illustration: Tate cohomology}
\label{sec:tate}

\noindent
In this final section, we demonstrate how one can utilize the
Gorenstein flat-cotorsion dimension to generalize a result of Hu and
Ding \cite{JHuNDn16}.

\begin{dfn}
  Let $M$ be an $R$-complex. A \emph{complete flat-cotorsion
    resolution} of $M$ consists of the following data: a
  semi-flat-cotorsion replacement $W$ of $M$, a totally acyclic
  complex of flat-cotorsion $R$-modules $T$, and a morphism
  $\mapdef{\tau}{T}{W}$ of $R$-complexes such that $\tau_i$ is an
  isomorphism for all $i\gg 0$.
\end{dfn}

\begin{prp}
  \label{prp:existence}
  Let $M$ be an $R$-complex and $n$ an integer. The following
  conditions are equivalent.
  \begin{eqc}
  \item $\Gfcd{M} \le n$.
  \item For every semi-flat-cotorsion replacement $W$ of $M$ there
    exists a complete flat-cotorsion resolution $\tau: T\to W$ of $M$
    such that $\tau_i$ is a split epimorphism for every $i\in\ZZ$ and
    $\tau_i=\id[W_i]$ holds for every $i\ge n$.
  \item There exists a complete flat-cotorsion resolution
    $\tau: T\to W$ of $M$ such that $\tau_i$ is an isomorphism for
    every $i\ge n$.
  \end{eqc}
\end{prp}

\begin{prf*}
  The proof is cyclic; the implication \proofofimp[]{ii}{iii} is
  trivial.

  \proofofimp{i}{ii} Let $W$ be a semi-flat-cotorsion replacement of
  $M$. By assumption, the module $\Co[n]{W}$ is Gorenstein
  flat-cotorsion, so it follows from \lemref{enough} that there is an
  acyclic complex concentrated in degrees $n$ and below,
  \begin{equation*}
    T' \deq 0 \lra \Co[n]{W} \lra T_{n-1} \lra T_{n-2} \lra \cdots\:,
  \end{equation*}
  where the modules $T_i$ are flat-cotorsion and such that
  $\Hom{T'}{V}$ is acyclic for every flat-cotorsion $R$-module
  $V$. Let $T$ be the complex obtained by splicing together
  $\Tha{n-1}{T'}$ and $\Thb{n}{W}$ at $\Co[n]{W}$. Now the identity on
  $\Co[n]{W}$ lifts to a morphism $\Tha{n-1}{T} \to \Tha{n-1}{W}$. The
  induced degreewise split surjective homomorphism
  $\Tha{n-1}{T} \oplus \Shift[-1]{\Cone{1^{\Tha{n-1}{W}}}} \to \Tha{n-1}{W}$
  together with the identity on $\Thb{n}{T} = \Thb{n}{W}$ is the
  desired morphism $\tau$.

  \proofofimp{iii}{i} As $\tau_i$ is an isomorphism for $i \ge n$ one
  has $\Co[n]{W} \is \Co[n]{T}$, and the latter module is Gorenstein
  flat-cotorsion.
\end{prf*}

\begin{lem}
  \label{lem:unique}
  Let $M$ be an $R$-complex of finite Gorenstein flat-cotorsion
  dimension. If $\mapdef{\tau}{T}{W}$ and $\mapdef{\tau'}{T'}{W'}$ are
  complete flat-cotorsion resolutions of $M$, then there is a homotopy
  equivalence $T \to T'$.
\end{lem}

\begin{prf*}
  For an $R$-module $X$ the disk complex $0 \to X\xra{=} X\to 0$
  concentrated in degrees $i+1$ and $i$ is denoted $\dc{i}{X}$.  Fix
  $i > \Gfcd{M}$; as $\Gfcd{M} \ge \sup{M}$ holds, one has
  $\Co[i+1]{W} \is \Cy[i]{W}$ and $\Co[i+1]{W'} \is \Cy[i]{W'}$. By \lemref{Schanuel} there exist
  flat-cotorsion modules $V$ and $V'$ with
  $\Cy[i]{W} \oplus V \is \Cy[i]{W'} \oplus V'$.  Set
  \begin{equation*}
    \widetilde{W} = W \oplus \dc{i}{V}\,, \
    \widetilde{W}' = W' \oplus \dc{i}{V'}\,, \
    \widetilde{T} = T \oplus \dc{i}{V}\,, \text{ and }
    \widetilde{T}' = T' \oplus \dc{i}{V'}\:.
  \end{equation*}
  The induced morphisms
  $\mapdef{\tilde\tau}{\widetilde{T}}{\widetilde{W}}$ and
  $\mapdef{\tilde\tau'}{\widetilde{T}'}{\widetilde{W}'}$ are complete
  flat-cotorsion resolutions of $M$. The isomorphism
  $\Cy[i]{\widetilde{W}} \is \Cy[i]{\widetilde{W'}}$ lifts by
  \lemcite[3.1 and Prop.~3.3]{CETa} to a homotopy equivalence
  $\widetilde{T} \to \widetilde{T}'$. Evidently, there are homotopy
  equivalences $T \to \widetilde{T}$ and $\widetilde{T}' \to T'$.
\end{prf*}

In view of \lemref{unique} and \prpref{compare} one can now define a
version of Tate cohomology.

\begin{dfn}
  \label{dfn:Tate}
  Let $M$ be an $R$-complex of finite Gorenstein flat-cotorsion
  dimension and $\mapdef{\tau}{T}{W}$ a complete flat-cotorsion
  resolution of $M$. For an $R$-complex $N$ and $i\in\ZZ$ set
  \begin{equation*}
    \Text{i}{M}{N} = \H[-i]{\Hom{T}{N}}\:.
  \end{equation*}
\end{dfn}

For a complex of finite Gorenstein flat dimension, every complete
flat-cotorsion resolution is per the next lemma a Tate flat resolution
in the sense of \cite{LLn13,CKL-17}.

\begin{lem}
  \label{lem:Ftac}
  Let $M$ be an $R$-complex of finite Gorenstein flat dimension. In
  every complete flat-cotorsion resolution $T \to W$ of $M$, the
  complex $T$ is \Ftac.
\end{lem}

\begin{prf*}
  Let $T \to W$ be a complete flat-cotorsion resolution of $M$. As $W$
  is a semi-flat replacement of $M$, the modules
  $\Co[i]{T} \is \Co[i]{W}$ are Gorenstein flat for $i \gg 0$. Thus
  $\Co[i]{T}$ is a Gorenstein flat-cotorsion module of finite
  Gorenstein flat dimension for every $i\in\ZZ$. As the dimensions
  agree by \thmref{Gfd-Gfcd}, it follows that $\Co[i]{T}$ is
  Gorenstein flat for every $i\in\ZZ$, whence $T$ is \Ftac, see
  \lemcite[2.3]{CFH-06}.
\end{prf*}

\begin{rmk}
  \label{rmk:HuDing}
  Let $M$ be an $R$-complex. By \fctref{gil} there is an exact
  sequence of $R$-complexes $0 \to M \to C \to F' \to 0$ with $C$
  semi-cotorsion and $F'$ acyclic and semi-flat and another exact
  sequence $0 \to C' \to F \to C \to 0$ with $F$ semi-flat and $C'$
  semi-cotorsion and acyclic. The diagram $F \qra C \qla M$ is a
  semi-flat-cotorsion replacement of $M$.  Assume that $M$ has finite
  Gorenstein flat dimension. By \thmref{Gfd-Gfcd} and
  \prpref{existence} there is a complete flat-cotorsion resolution
  $\mapdef{\tau}{T}{W}$ with $\tau_i$ a split epimorphism for all
  $i\in\ZZ$. By \prpref{compare} there is a homotopy equivalence
  $W \to F$, so by \lemref{Ftac} the diagram
  \begin{equation*}
    T \xra{\tau} W \qra C \qla M
  \end{equation*}
  is a complete flat resolution in the sense of \cite{JHuNDn16}.  It
  follows that for $R$, $M$, and $N$ as in \thmcite[5.5]{JHuNDn16},
  our \dfnref{Tate} agrees with the definition in \cite{JHuNDn16}.
\end{rmk}

In view of \thmref{Gfd-Gfcd} and \rmkref{HuDing} the next result
generalizes parts of \thmcite[1.5]{JHuNDn16}.

\begin{thm}
  \label{thm:vanishing}
  Let $M$ be an $R$-complex. If $M$ has finite Gorenstein
  flat-cotorsion dimension, then the following statements are
  equivalent.
  \begin{eqc}
  \item $\fd{M} = \Gfcd{M}$ holds.
  \item $\Text{m}{M}{N}=0$ holds for all $m\in\ZZ$ and every
    $R$-complex $N$.
  \item $\Text{m}{M}{C}=0$ holds for some $m\in\ZZ$ and every
    cotorsion $R$-module $C$.
  \end{eqc}
\end{thm}

\begin{prf*}
  The proof is cyclic; the implication \proofofimp[]{ii}{iii} is
  clear.

  \proofofimp{i}{ii} By \lemref{fdlemma} there is a
  semi-flat-cotorsion replacement $W$ of $M$ with $W_i = 0$ for all
  $i \gg 0$. It follows that $0\to W$ is a complete flat-cotorsion
  resolution of $M$.

  \proofofimp{iii}{i} Since $M$ has finite Gorenstein flat-cotorsion
  dimension, there exists by \prpref{existence} a complete
  flat-cotorsion resolution $\mapdef{\tau}{T}{W}$ of $M$ with
  $\tau_i=\id[F_i]$ for all $i\ge g$, where $g$ is some integer. The
  module $\Co[m]{T}$ is cotorsion, so
  $\Text{m}{M}{\Co[m]{T}} = \H[-m]{\Hom{T}{\Co[m]{T}}}=0$ holds. This
  means that application of $\Hom{-}{\Co[m]{T}}$ leaves the sequence
  \begin{equation*}
    0 \lra \Co[m]{T} \lra T_{m-1} \lra \Co[m-1]{T} \lra 0
  \end{equation*}
  exact, and it follows that it splits. In
  particular, $\Co[m-1]{T}$ is a flat module and hence so is
  $\Co[i]{T}$ for every $i \ge m-1$. For $i \gg 0$ the module
  $\Co[i]{W} = \Co[i]{T}$ is flat, so $M$ has finite flat
  dimension. Now invoke \thmref{fd-Gfcd}.
\end{prf*}

\section*{Acknowledgment}

\noindent
We thank the anonymous referee for numerous pertinent comments that
helped improve the exposition. In particular, we are thankful for
suggestions that shaved half a page off our original proof of
\prpref{ses}(a).


\def\soft#1{\leavevmode\setbox0=\hbox{h}\dimen7=\ht0\advance \dimen7
  by-1ex\relax\if t#1\relax\rlap{\raise.6\dimen7
  \hbox{\kern.3ex\char'47}}#1\relax\else\if T#1\relax
  \rlap{\raise.5\dimen7\hbox{\kern1.3ex\char'47}}#1\relax \else\if
  d#1\relax\rlap{\raise.5\dimen7\hbox{\kern.9ex \char'47}}#1\relax\else\if
  D#1\relax\rlap{\raise.5\dimen7 \hbox{\kern1.4ex\char'47}}#1\relax\else\if
  l#1\relax \rlap{\raise.5\dimen7\hbox{\kern.4ex\char'47}}#1\relax \else\if
  L#1\relax\rlap{\raise.5\dimen7\hbox{\kern.7ex
  \char'47}}#1\relax\else\message{accent \string\soft \space #1 not
  defined!}#1\relax\fi\fi\fi\fi\fi\fi}
  \providecommand{\MR}[1]{\mbox{\href{http://www.ams.org/mathscinet-getitem?mr=#1}{#1}}}
  \renewcommand{\MR}[1]{\mbox{\href{http://www.ams.org/mathscinet-getitem?mr=#1}{#1}}}
  \providecommand{\arxiv}[2][AC]{\mbox{\href{http://arxiv.org/abs/#2}{\sf
  arXiv:#2 [math.#1]}}} \def\cprime{$'$}
\providecommand{\bysame}{\leavevmode\hbox to3em{\hrulefill}\thinspace}
\providecommand{\MR}{\relax\ifhmode\unskip\space\fi MR }
\providecommand{\MRhref}[2]{%
  \href{http://www.ams.org/mathscinet-getitem?mr=#1}{#2}
}
\providecommand{\href}[2]{#2}

\end{document}